\documentclass[a4paper,10pt]{article}
\newcommand{\colorprint}[1]{}
\usepackage[utf8]{inputenc}
\usepackage{color}
\usepackage{tikz}
\usetikzlibrary{shapes}
\usepackage{float}
\usepackage{subcaption}
\usepackage{amsmath}
\usepackage{amsfonts}
\usepackage{amssymb}
\usepackage{amsthm}
\usepackage{mathrsfs}
\usepackage{amsthm}
\usepackage{subcaption}
\usepackage{graphics}
\usepackage[colorlinks,urlcolor=magenta,linkcolor=blue,citecolor=teal]{hyperref}
\theoremstyle{plain}
\newtheorem{theorem}{Theorem}

\newtheorem{lemma}[theorem]{Lemma}

\newtheorem{conjecture}[theorem]{Conjecture}
\newcommand{\dom}{\rightarrow}

\newcommand{\induce}[2]{#1\langle{}#2\rangle}

   \newcommand{\NoColor}[1]{}

\begin{document}
\bibliographystyle{plain}
\title{Every $(13k-6)$-strong tournament with minimum out-degree at least $28k-13$ is $k$-linked}

\author{J\o{}rgen Bang-Jensen and Kasper Skov Johansen\thanks{Department of Mathematics and Computer Science, University of Southern Denmark, email (jbj@imada.sdu.dk and kajoh14@student.sdu.dk). Research supported by the Independent Research Fund Denmark under grant number DFF 7014-00037B. This paper was written while the second author was a masters student at Department of Mathematics and Computer Science}}
\date{June 23, 2021}
\maketitle
\begin{abstract}
  A digraph $D$ is $k$-linked if it satisfies that  for every choice of disjoint sets $\{x_1,\ldots{},x_k\}$ and $\{y_1,\ldots{},y_k\}$ of vertices of $D$ there are vertex disjoint paths $P_1,\ldots{},P_k$ such that $P_i$ is an $(x_i,y_i)$-path. Confirming a conjecture by K\"uhn et al, Pokrovskiy proved in 2015 that every $452k$-strong  tournament is $k$-linked and asked for a better linear bound. Very recently Meng et al proved that every $(40k-31)$-strong tournament is $k$-linked. In this note  we use an important lemma from  their paper to give a short  proof that every $(13k-6)$-strong tournament of minimum out-degree at least $28k-13$ is $k$-linked.
\end{abstract}
\noindent{}{\bf Keywords:} tournament, linkage, connectivity.

\section{Introduction}
A digraph $D$ is {\bf strongly connected} if it has a directed path from $x$ to $y$ (an $
(x,y)$-path) for every ordered pair of distinct vertices $x,y$ and it is {\bf $\mathbf{k}$-strong} if it has at least $k+1$ vertices and remains strongly connected when we delete any set of at most $k-1$ vertices. A (di)graph $D$ is {\bf $\mathbf{k}$-linked} if it has  vertex disjoint paths $P_1,\ldots{},P_k$ such that $P_i$ is an $(x_i,y_i)$-path for every choice of disjoint sets $\{x_1,\ldots{},x_k\}$ and $\{y_1,\ldots{},y_k\}$ of vertices of $D$.
Thomas and Wollan proved that for every natural number $k$ every $10k$-connected undirected graph is $k$-linked \cite{thomasEJC26}. Hence there exists a smallest function $f(k)$ such that every $f(k)$-connected undirected graph is $k$-linked.

Thomassen \cite{thomassenC11} showed that for general digraphs there is no integer $\ell$ so that every $\ell$-strong digraph is 2-linked, thus ruling out the existence of a similar function $f(k)$ for directed $k$-linkage. Hence it is natural  to focus on special classes of digraphs and one such important class of digraphs is tournaments, that is, digraphs in which there is exactly one arc between every pair of distinct vertices. Thomassen \cite{thomassen1984} was the first to prove that for tournaments there exists a function $h(k)$ such that every $h(k)$-strong tournament is $k$-linked. Bang-Jensen \cite{bangADM41} proved that $h(2)=5$ but already the value of $h(3)$ is open. Thomassen's
function $h(k)$ was exponential in $k$. This was improved to a polynomial in $k$ by K\"uhn, Lapinskas, Osthus, and Patel  in \cite{kuehnPLMS109} and they conjectured that a linear function would suffice. Pokrovskiy confirmed this  in \cite{pokrovskiyJCT115} by showing that every $452k$-strong tournament is $k$-linked. 

If $X=\{x_1,\ldots{},x_k\}$ and $Y=\{y_1,\ldots,y_k\}$ are disjoint subsets of vertices in a digraph $D$, we say that $X$ {\bf anchors} $Y$ in $D$ if for every permutation $\pi$ of $\{1,2,\ldots,k\}$ there exists a collection of disjoint paths $Q_1,\ldots{},Q_k$ in $D$ such that $Q_i$ is an $(x_i,y_{\pi(i)})$-path.
The proof in \cite{pokrovskiyJCT115} is based on the following lemma which shows that tournaments have a rich structure.

\begin{lemma}\cite{pokrovskiyJCT115}
  \label{lem:11k}
  For every narural number $p$ every tournament $T$ on at least $11p$ vertices contains disjoints sets $X,Y$ both of size $p$ such that $X$ anchors $Y$ in $T$ 
  \end{lemma}

  Very recently Meng, Rolek, Wang and Yu \cite{mengEJCta} used the following  improved version of Lemma \ref{lem:11k} together with some nice ideas to show that every $(40k-31)$-strong tournament is $k$-linked.

  \begin{lemma}\cite{mengEJCta}
    \label{lem:9k}
    Let $p$ be a natural number and let $T$ be a tournament on $n\geq 9p-6$ vertices. Then $T$ contains a pair of disjoint vertex sets $X,Y$, each of size $p$ such that $X$ anchors $Y$ in $T$.
    \end{lemma}

  In this paper we use Lemma \ref{lem:9k} to obtain the following better bound.

\begin{theorem}
  \label{thm:main}
  For every natural number $k$ every $(13k-6)$-strong tournament of minimum out-degree at least $28k-13$ is $k$-linked.
\end{theorem}

The proof in \cite{mengEJCta} that every $(40k-31)$-strong tournament is $k$-linked  uses an approach similar to that used by Pokrovskiy in \cite{pokrovskiyJCT115} but is based on Lemma \ref{lem:9k} instead of Lemma \ref{lem:11k} and some matching arguments. Our proof also uses Lemma \ref{lem:9k} but we avoid the use of matching arguments.

Notation not already introduced or given below  is consistent with \cite{bang2018,bang2009}.
Let $D=(V,A)$ be a digraph and let $X\subseteq V$, then we denote by $D\langle{}X\rangle$ the subdigraph of $D$ induced by $X$. We denote by $N_D^+(v), N_D^-(v)$ the set of out-neighbours, respectively in-neighbours of $v$ in $D$, that is, the set of vertices $y\neq v$ such that $v\rightarrow y$, respectively $y\rightarrow v$ is an arc of $D$. We call 
$d^+_D(v)=|N_D^+(v)|$,  $d^-_D(v)=|N_D^-(v)|$ the {\bf out-degree}, respectively the {\bf in-degree} of the vertex $v$ in $D$. If $X$ is a subset of the vertices of $D$ then we use the notation $D-X$ for the digraph $\induce{D}{V(D)\setminus X}$.

\section{Proof of Theorem \ref{thm:main}}

Let $T_0 = T - (X_0 \cup Y_0)$. Pick  a vertex $u_1\in V(T_0)$ of minimum out-degree in $T_0$, that is $d^+_{T_0}(u_1)\leq d^+_{T_0}(v)$ for all $v\in V(T_0)$. Then we find a vertex $v_1$ with minimum  out-degree in $T\langle N^+_{T_0}(u_1)\rangle$. Now let $A_1 = N^+_{T_0}(u_1) \cap N^+_{T_0}(v_1)$ and set $D_1 = \lbrace u_1, v_1 \rbrace$. Also, let $T_1=T_0 - D_1$.\\
\newline
For $i = 2,...,9k-6$ in that order we pick $u_i, v_i \in T_{i-1}$ such that $u_i$ has minimum out-degree in $T_{i-1}$ and $v_i$ is a vertex of minimum out-degree in $T\langle N^+_{T_{i-1}}(u_i) \rangle$. Let $A_i = N^+_{T_{i-1}}(u_i) \cap N^+_{T_{i-1}}(v_i)$, $D_i = \lbrace u_i, v_i \rbrace$ and $T_i = T_{i-1} - D_i$. Finally set $U = \lbrace u_1,...,u_{9k-6} \rbrace$ and $V = \lbrace v_1,...,v_{9k-6} \rbrace$.\\

The tournament $T\langle N^+_{T_{i-1}}(u_i) \rangle$ has $d^+_{T_{i-1}}(u_i)$ vertices so it has a vertex $w$ whose  out-degree in $T\langle N^+_{T_{i-1}}(u_i) \rangle$ 
is at most $(d^+_{T_{i-1}}(u_i)-1)/2$. Because $v_i$ was picked to have smallest out-degree in $T\langle N^+_{T_{i-1}}(u_i)\rangle$ it follows that

\begin{equation}\label{1st}
\vert A_i \vert = \vert N^+_{T_{i-1}}(u_i) \cap N^+_{T_{i-1}}(v_i) \vert = d^+_{T\langle N^+_{T_{i-1}}(u_i)\rangle}(v_i) \leq (d^+_{T_{i-1}}(u_i)-1)/2.
\end{equation}\\
By (\ref{1st}) and the fact that $u_i$ has minimum out-degree in $T_{i-1}$ we get

\begin{equation}\label{2nd}
\vert A_i \vert \leq (d^+_{T_{i-1}}(u_i)-1)/2 \leq (d^+_{T_{i-1}}(z)-1)/2, \quad \forall z\in A_i
\end{equation}\\
Then since $V(T_{i-1}) \subseteq V(T)$ we get from (\ref{2nd}) that

\begin{equation}\label{3rd}
\vert A_i \vert \leq (d^+_{T_{i-1}}(z)-1)/2 \leq (d^+_T(z)-1)/2, \quad \forall z\in A_i.
\end{equation}\\
Since $V$ contains $9k-6$ vertices we know, by Lemma \ref{lem:9k}, that some disjoint subsets $V'=\lbrace v'_1,...,v'_k \rbrace$ and $V''= \lbrace v''_1,...,v''_k \rbrace$ exist in $V$ such that $V'$ anchors $V''$ within $T\langle V \rangle$. Let $V^* = V-(V'\cup V'')$. We also collect the $k$ vertices $u'_i$ from $U$, for $i=1,...,k$, where $v'_i$ was picked in the tournament induced on the out-neighbors of $u'_i$. Let $U' = \lbrace u'_1,..., u'_k \rbrace$ and $U^* = U - U'$. This guarantees that $u'_i\dom v'_i$ is an arc and that $D_{j_i} = \lbrace u'_i, v'_i \rbrace$ for some $j_i\in \lbrace 1,...,9k-6 \rbrace$. Then we let $A'_i = A_{j_i}$ and we get from (\ref{3rd})  that
\begin{equation}\label{4th}
\vert A'_i \vert \leq (d^+_T(z)-1)/2, \quad \forall z\in A'_i.
\end{equation}\\
Next we consider the subtournament  $T' = T-(X_0 \cup Y_0 \cup U' \cup V)$. Because $\vert X_0 \vert + \vert Y_0 \vert + \vert U' \vert + \vert V \vert = 12k-6$ it follows that each  $x_i \in X_0$ has at least $16k-7$ out-neighbors in $V(T')$ (notice that $x_i$ is not in $V(T')$ itself). We can therefore pick a set $X' = \lbrace x'_1,...,x'_k \rbrace$ of $k$ distinct vertices   of $V(T')$ such that $x_i \rightarrow x'_i$ is an arc of $T$.\\
\newline
Now we are going to construct $k$ disjoint paths $P_1,\ldots{},P_k$   such that $P_i$ is an $(x_i,v'_i)$-path and these paths are disjoint.

Let $I =  \lbrace 1,2,...,k \rbrace$ and let $\hat{I} \subseteq I $ be those indexes for which $x'_i$ has an arc to at least one of $u'_i,v'_i$. 
If $i \in \hat{I}$ then $T$ contains either the path $x_i\dom x'_i\dom v'_i$ or the path $x_i\dom x'_i \dom u'_i\dom v'_i$.

For $i \in I - \hat{I}$ we have that $x'_i$ is in $A'_i$, and by (\ref{4th}) it follows that 
\begin{equation}\label{5th}
\vert A'_i \vert \leq (d^+_T(x'_i)-1)/2.
\end{equation}
Consider now the out-neighbors of $x'_i$ in $T$. We have
\begin{equation}\label{6th}
\vert X_0 \vert + \vert Y_0 \vert + \vert U' \vert + \vert V \vert + \vert X' \vert + k - 1 = 14k-7 = (28k-14)/2 \leq (d^+_T(x'_i) - 1)/2
\end{equation}
Then by adding (\ref{5th}) and (\ref{6th}) we get 
\begin{equation}
\vert A'_i \vert + \vert X_0 \vert + \vert Y_0 \vert + \vert U' \vert + \vert V \vert + \vert X' \vert + k - 1 \leq d^+_T(x'_i) - 1.
\end{equation}
That is, 
\begin{equation}
\vert A'_i \vert + \vert X_0 \vert + \vert Y_0 \vert + \vert U' \vert + \vert V \vert + \vert X' \vert + k \leq d^+_T(x'_i).
\end{equation}
In other words, we see that $x'_i$ has at least $k$ out-neighbors outside of $(A'_i \cup X_0 \cup Y_0 \cup U' \cup V \cup X')$. We can therefore find a set $X'' = \lbrace x''_i \vert i \in I-\hat{I} \rbrace$ of distinct vetices   with $X'' \subseteq V(T')-X'$ such that $x'_i \rightarrow x''_i$ is an arc in $T$ and such that at least one of the arcs $x''_i \rightarrow u'_i$ and $x''_i \rightarrow v'_i$ exists in $T$. The latter follows by the fact that for all 
$i\in I-\hat{I}$ the vertex $x''_i$ is not in $A'_i$. Therefore we can find, for all $i=1,...,k$,
 some path from $x_i$ to $v'_i$ that has one of the following forms:
\begin{itemize}
\item $x_i \rightarrow x'_i \rightarrow v'_i$
\item $x_i \rightarrow x'_i \rightarrow u'_i \rightarrow v'_i$
\item $x_i \rightarrow x'_i \rightarrow x''_i \rightarrow v'_i$
\item $x_i \rightarrow x'_i \rightarrow x''_i \rightarrow u'_i \rightarrow v'_i$
\end{itemize}

Call such paths $P_1,...,P_k$ and note that, by construction, they are disjoint.
Then each $P_i$ has at most two vertices in $V(T')$, namely $x'_i$ and $x''_i$. All vertices of the
paths $P_i$ are contained in $(X_0 \cup U' \cup V' \cup X' \cup X'')$. Furthermore, $\vert X_0 \cup U' \cup V' \cup X' \cup X'' \vert \leq \vert X_0 \cup U' \cup V' \cup V^* \cup X' \cup X'' \vert \leq 12k-6$. Hence, since $T$ is $(13k-6)$-strong it follows that $T^* = T-(X_0 \cup U' \cup V' \cup V^* \cup X' \cup X'')$ is $k$-strong and contains both $Y_0$ and $V''$. By Menger's Theorem there exist $k$ disjoint paths from $V''$ to $Y_0$ in $T^*$. Name these $R_1,...,R_k$ so that $R_i$ ends in $y_i$ for $i\in \{1,2,\dots{},k\}$. Now, neither of the paths $R_j$ intersect some path $P_i$ because the vertices of each $P_i$ are  completely contained in $(X_0 \cup U' \cup V' \cup X' \cup X'')$. Suppose now we pair up $P_i$ and $R_i$ for each $i\in \lbrace 1,...,k \rbrace$. Let $h(P_i)$ denote the terminal (last) vertex of $P_i$ and let $t(R_j)$ denote the initial vertex of $R_j$. Then $(h(P_1),...,h(P_k))$ is a permutation of $V'$ and $(t(R_1),..., t(R_k))$ is a permutation of $V''$. By the choice of $V'$ and $V''$ we can find a linkage $\lbrace M_1,...,M_k \rbrace$ from $\lbrace h(P_1),...,h(P_k) \rbrace$ to $\lbrace t(R_1),..., t(R_k) \rbrace$ whose inner vertices are completely within $V^*$. But then since $P_i$ and $R_i$ are disjoint from $V^*$ it is easy to see that $Q_1, Q_2, ..., Q_k$ is a linkage from $X_0$ to $Y_0$ when we set $Q_i = P_i M_i R_i$ for each $i \in \{1,2,\dots{},k\}$.

\qed
\section{Remarks}
Pokrovskiy gave an infinite family of $(2k-2)$-strong tournaments which  have arbitrary high minimum in-degree and out-degree but are not $k$-linked and made following conjecture.

\begin{conjecture}\cite{pokrovskiyJCT115}
  \label{conj:2klinked}
  For every natural number $k$ there exists an integer $d(k)$ such that every $2k$-strong tournament
  with minimum in-degree and minimum out-degree at least $d(k)$ is $k$-linked.
\end{conjecture}

The following partial result was obtained by Gir{\~{a}}o  and Snyder

\begin{theorem}\cite{giraoJCT139}
  \label{thm:girao}
  For every natural number $k$ there exists an integer $d(k)$ such that every $4k$-strong tournament
  with minimum out-degree at least $d(k)$ is $k$-linked.
\end{theorem}

In an unpublished paper \cite{giraoarXiv1912.00710} Gir{\~{a}}o, Popielarz and Snyder proved Conjecture \ref{conj:2klinked} up to an additive constant of 1.

\begin{theorem}\cite{giraoarXiv1912.00710}
  \label{thm:2k+1}
  For every natural number $k$ there exists an integer $d(k)$ such that every $(2k+1)$-strong tournament with minimum out-degree at least some polynomial in k is k-linked.
\end{theorem}

The function $d(k)$ used in Theorem \ref{thm:2k+1} is $d(k)=Ck^{31}$ for some appropriate constant $C$ and hence Theorem \ref{thm:2k+1} does not imply Theorem \ref{thm:main}.

The following result shows that the condition on high minimum out-degree cannot be removed in Conjecture \ref{conj:2klinked}.

\begin{theorem}\cite{giraoarXiv1912.00710}
  For every natural number $k$ there exists infinitely many  $(5k-1)$-strong tournaments which are  not $2k$-linked.
  \end{theorem}

  Gir{\~{a}}o, Popielarz and Snyder also showed that $2k$ cannot be replaced by $2k-1$ in Conjecture \ref{conj:2klinked}.
  
  \begin{theorem}\cite{giraoarXiv1912.00710}
    For all integers $k\geq 2$ and $m\geq 2k$ there exists a $(2k-1)$-strong tournament $T$  whose minimum in-degree and out-degree is at least $m$ such that $T$ is not $k$-linked.
    \end{theorem}

    A digraph is {\bf semicomplete} if it has no pair of non-adjacent vertices. Thus the class of tournaments forms a subclass of the class of semicomplete digraphs. It was pointed out in \cite{bangTchapter}  that Pokrovskiy's proof in \cite{pokrovskiyJCT115} also holds for semicomplete digraphs, implying that every $452k$-strong semicomplete digraph is $k$-linked. The proof of Theorem \ref{thm:main} does not hold for semicomplete digraphs because we cannot guarantee a vertex of out-degree at most $(d^+_{T_{i-1}}(u_i)-1)/2$ in $T\langle N^+_{T_{i-1}}(u_i) \rangle$ (just above (\ref{1st})) as the digraph $T\langle N^+_{T_{i-1}}(u_i) \rangle$ could be complete. It can also be checked that the proof in \cite{mengEJCta} does not hold for semicomplete digraphs.

    It is known \cite{guoDAM79b} that every $(3k-2)$-strong semicomplete digraph contains a spanning $k$-strong tournament. From this and Theorem \ref{thm:main} it follows that every
    $(84k-41)$-strong semicomplete digraph is $k$-linked but we believe this bound is far from best possible. We also believe that Conjecture \ref{conj:2klinked}
   holds for semicomplete digraphs.

\end{document}